\documentclass[11pt]{amsart}
\usepackage{amssymb}
\usepackage{amsthm}
\usepackage{enumerate}
\usepackage{cite}

\renewcommand{\phi}{\varphi}
\renewcommand{\Re}{\mathop{\mathrm{Re}}\nolimits}
\renewcommand{\le}{\leqslant}
\renewcommand{\ge}{\geqslant}

\newcommand{\eps}{\varepsilon}
\newcommand{\T}{\mathbb T}
\newcommand{\N}{\mathbb N}
\newcommand{\Z}{\mathbb Z}

\newcommand{\C}{\mathbb C}

\newcommand{\sgn}{\mathop{\mathrm{sgn}}\nolimits}
\newcommand{\rank}{\mathop{\mathrm{rank}}\nolimits}
\newcommand{\ARG}{\mathop{\mathrm{arg}}\nolimits}

\newtheorem{Thm}{Theorem}
\newtheorem{Lem}[Thm]{Lemma}
\newtheorem{Prop}[Thm]{Proposition}

\newtheorem*{Def}{Definition}

\def\beginpf{\medskip\noindent {\bf Proof.} \;}

\title[The past and future wave operators.]
{The past and future wave operators \\ on the singular spectrum.}

\author{R.\,V.\,Bessonov}

\address{St.Petersburg Department of Steklov Mathematical Institute RAS,
Fontanka 27, St.Petersburg, 191023, Russia}

\address{Chebyshev Laboratory, St.-Petersburg State University,
14th Line, 29b, Saint-Petersburg, 199178 Russia}

\email{bessonov@pdmi.ras.ru}


\thanks{The author is partially supported by the RFBR grant 11-01-00584-a, by the V.~A.~Rokhlin grant and by the Chebyshev Laboratory
(Department of Mathematics and Mechanics, St.-Petersburg State University)
under RF government grant 11.G34.31.0026}

\begin{document}

\maketitle
\sloppy
\section{Introduction}\label{introd}

Given a pair of Hilbert spaces $H_1, H_2$, unitary operators $U_1 :H_1 \to H_1$, $U_2: H_2 \to H_2$, and a bounded operator $A : H_1 \to H_2$, define  the future and past wave operators $W_{+}$, $W_{-}$ as the limits of the sequence $\{U_{2}^{-n}AU_{1}^{n}\}_{n \in \Z}$  as $n \to +\infty$, $n \to -\infty$, respectively.  

There are a number of different meanings of the above term ``limit''. In the classical scattering theory it is proved that the limits $\lim_{n \to \pm \infty} U_{2}^{-n}AU_{1}^{n}$ exist in the strong operator topology, under the assumptions that the spectral measures of $U_1$, $U_2$ are absolutely continuous with respect to the Lebesgue measure and the commutator $AU_1 - U_2A$ is of trace class\cite{CST}.  

Much less is known about similar results for unitary operators $U_1$,  $U_2$ having singular spectrum. Simple examples in a one-dimensional space show that the limits $\lim_{n \to \pm \infty} U_{2}^{-n}AU_{1}^{n}$ may not exist. A natural way to define the wave operators in this situation is to consider the limits of \textit{averages} of the sequence $U_{2}^{-n}AU_{1}^{n}$ by using some summation (averaging) method. Important examples of summation methods which will be used here are the Cesaro summation method taking a sequence to the sequence of its arithmetical means and the Abel-Poisson summation method.

\medskip 

Now we give a precise definition of summation methods and  averaged wave operators. Let $(E, \prec)$ be a totally ordered set, $\Z_+$ the set of nonnegative integers. A family of nonnegative numbers $\{p_{\alpha,n}\}_{\alpha \in E, n \in \Z_+}$ determines an s-regular summation (averaging) method if the following conditions are fulfilled: 

$1)$ $\sum_{n} p_{\alpha,n} = 1$ for every $\alpha \in E$,

$2)$ $\lim_\alpha \sum_{n} |p_{\alpha,n} - p_{\alpha,n+1}| = 0$,

$3)$ $\lim_\alpha p_{\alpha,n} = 0$ for every $n \in \Z_+$.

\noindent For a sequence $\{x_n\}$ of elements of a Banach space $X$, which is bounded in norm, define its averages by $x(\alpha) = \sum_n p_{\alpha,n} x_n$, $\alpha \in E$. If the limit $x = \lim_\alpha x(\alpha)$ exists in some topology $\sigma$ on $X$, we will say that the averages of $\{x_n\}$ converge to $x$ in $\sigma$. Usually we will work with the space $X = B(H)$ of all bounded operators on a Hilbert space $H$ equipped by the weak operator topology $\sigma$. All s-regular methods possess the following natural properties: 
\begin{itemize}
\item[$a)$] Regularity: averages of every convergent sequence converge to its limit. 
\item[$b)$] Linearity: if averages of sequences $\{x_n\}$, $\{y_n\}$ converge to $x$, $y$, then averages of $\{x_n + cy_n\}$ converge to $x + cy$ for any complex \mbox{number $c$;}
\item[$c)$] Stability: averages of sequences $\{x_n\}$ and $\{x_{n+1}\}$
 do or do not converge simultaneously, and if they converge, then the corresponding limits coincide.   
\end{itemize}
Moreover, every s-regular summation method averages to zero any nonconstant unimodular geometric progression of complex numbers:
\begin{itemize}
\item[$d)$] $\lim_\alpha \sum_{n = 0}^{\infty} p_{\alpha,n} z^n = 0$ if $z \in \C \setminus \{1\}$, $|z| = 1$.
\end{itemize}
The well-known examples of s-regular methods are the Cesaro summation method (with $E = (\N, \le)$ and $p_{\alpha,n} = \frac{1}{\alpha}$ whenever $n < \alpha$, $p_{\alpha,n} = 0$ otherwise), and the Abel-Poisson summation method (with $E = ((0,1), \le)$ and  $p_{\alpha,n} = (1-\alpha)\alpha^{n}$). 

\medskip

Fix some s-regular summation method. In what follows the term ``averaging'' will always mean the use of the fixed summation method. 

\medskip

\begin{Def}
Let $H_{1,2}$, $U_{1,2}$, $A$ be as above, denote by $W_{+}(\alpha)$, $W_{-}(\alpha)$ the averages of the sequences $\{U_{2}^{-n}AU_{1}^{n}\}_{n \in \Z_{+}}$ and $\{U_{2}^{n}AU_{1}^{-n}\}_{n \in \Z_{+}}$:
\begin{equation}
\label{wo}
W_{+}(\alpha) = \sum_{n = 0}^{+\infty} p_{\alpha,n} U_{2}^{-n}AU_{1}^{n}
\mbox{\quad and \quad}
W_{-}(\alpha) = \sum_{n = 0}^{+\infty} p_{\alpha,n} U_{2}^{n}AU_{1}^{-n}
.
\end{equation}
The future weak averaged wave operator $W_+$ is the limit $\lim_{\alpha}W_{+}(\alpha)$ in the weak operator topology if the limit exists. Similarly, the past weak averaged wave operator $W_-$ is the limit  $\lim_{\alpha}W_{-}(\alpha)$.  
\end{Def}         
We borrow the term ``wave operator'' from the classical scattering theory, where the case of absolutely continuous spectrum is studied.\\

If $AU_1 = U_2A$ we have $U_{2}^{-n}AU_{1}^{n} = A$, therefore operators $W_{\pm}$ obviously exist and equal $A$.  It seems natural to examine the problem of existence of $W_{\pm}$ for operators $U_1$, $U_2$, $A$ with a ``small'' commutator $AU_1 - U_2A$. A general conjecture can be formulated as follows: If $AU_1 - U_2A$ is a finite-rank operator, then the future and past averaged wave operators exist for every s-regular summation method. For the case of rank-one commutators this conjecture was essentially proven in \cite{vk}, see also \cite{vk-zns}. We study the case of rank-two commutators. Since the problem is solved for operators with absolutely continuous spectrum, one can restrict the consideration to singular unitary operators $U_1$, $U_2$. Our main result is the following: 
\begin{Thm}
\label{thm1}
Let $U_1: H_1 \to H_1$, $U_2: H_2 \to H_2$ be singular unitary operators, and let $A: H_1 \to H_2$ be a bounded operator such that $\rank(AU_1-U_2A)  \leqslant 2$. Then the weak limit $\lim_{\alpha}(W_{+}(\alpha) - W_{-}(\alpha))$ exists and equals zero. In particular, the weak averaged wave operators $W_{\pm}$ do or do not exist simultaneously. If $W_{\pm}$ exist, they coincide.  
\end{Thm} 

In this paper we follow the approach developed in \cite{vk-zns}. 
The general problem of existence of the limits $\lim_\alpha W_{\pm}(\alpha)$ and $\lim_\alpha (W_{+}(\alpha) - W_{-}(\alpha))$ for rank-two commutators reduces (see \cite[Theorem 7.2]{vk-zns}) to the following particular case:  
\begin{itemize}
\item[1)] $H_1 = H_2 = L^2(\mu)$, where $\mu$ is a Borel singular measure  on the unit circle $\T$ of the complex plane $\C$. The measure $\mu$ is free of point masses. 
\item[2)] $U_1 = U_2 = U$ is the operator of multiplication by the independent variable on $L^2(\mu)$,  
\item[3)] $AU - UA = (\cdot, \phi)1 - (\cdot, 1)\phi$ \; for some real-valued function $\phi \in L^2(\mu)$.
\end{itemize}
In this situation the problem of existence of wave operators $W_{\pm}$ can be restated in terms of functions $\phi$ from item $3)$. In \cite{VKnew} one can find a sufficient condition on the function $\phi$ guaranteeing the existence of $W_{\pm}$. Namely, if $\mu$  is the Clark measure $\sigma_1$ for an inner function $\theta$ and $\phi$ coincides $\mu$-almost everywhere with a trace of some function $h \in K_{\theta}$ having a continuous trace on $\sigma_{\alpha}$, $\alpha \neq 1$, then the operator $W_{-}$ exists. The precise definitions will be given in Section \ref{bb}, where we discuss the related results. 

We show that it suffices to consider only functions $\phi$ that coincide $\sigma_1$-almost everywhere with a continuous function on $\T$.
\begin{Thm} \label{thm2}
Suppose that for every triple $\mu$, $U$, $A$ satisfying conditions $1)-3)$, where, moreover, 
\begin{itemize}
\item[4)] $\phi$ coincides $\mu$-almost everywhere with a continuous function on $\T$,
\end{itemize}
one of the limits $\lim_\alpha W_{+}(\alpha)$, $\lim_\alpha W_{-}(\alpha)$,  $\lim_\alpha (W_{+}(\alpha) - W_{-}(\alpha))$ exists in the weak  operator topology.  Then the same limit exists for any operators $A$, $U_1$, $U_2$ from Theorem \ref{thm1} in the general case of $\rank(AU_1-U_2A) = 2$.
\end{Thm}
However, not all continuous functions arise from commutators of type $1)-3)$. Denote by $C(\T)$ the set of all continuous functions on the unit circle $\T$.
\begin{Thm} \label{thm3}
There exist a function $\phi \in C(\T)$ and a singular measure $\mu$ without atomic masses such that the operator $(\cdot, \phi)1 - (\cdot, 1)\phi$ on $L^2(\mu)$ cannot be represented as a commutator $AU - UA$ for a bounded operator $A$.\\
\end{Thm}
From the proof of Theorem \ref{thm3} it follows some consequences concerning to a boundary behaviour of pseudocontinuable functions. We discuss them in Section \ref{bb}.

\section{Proof of Theorem \ref{thm1}.} \label{proofthm1} For the proof of Theorem \ref{thm1} we need some preliminary technical results. In what follows we will assume that for the operators $A,U$ conditions $1)$ and $2)$ from Section \ref{introd} are fulfilled.  Define operators $W_{\pm}(\alpha)$ by formula \eqref{wo}. 

\begin{Lem} \label{diff} Let $K$ be the commutator $AU - UA$. We have   $$W_{+}(\alpha)U - UW_{-}(\alpha) = \sum_{n = 0}^{\infty} p_{\alpha, n} \sum_{l = -n}^{n} U^{l}KU^{-l}.$$ 
\end{Lem}
\beginpf We have
\begin{equation*}
\begin{split}
W_{+}(\alpha)U - UW_{-}(\alpha) & = \sum_{n = 0}^{\infty} p_{\alpha,n}U^{-n}AU^{n+1} - \sum_{n = 0}^{\infty} p_{\alpha, n}U^{n+1}AU^{-n}  \\&
= \sum_{n = 0}^{\infty} p_{\alpha,n}(U^{-n}AU^{n+1} - U^{n+1}AU^{-n})  \\&
= \sum_{n = 0}^{\infty} p_{\alpha,n} \sum_{-n}^{n} (U^{l}AU^{-l+1} - U^{l+1}AU^{-l})  \\& =
\sum_{n = 0}^{\infty} p_{\alpha, n} \sum_{l = -n}^{n} U^{l}KU^{-l}.
\end{split}
\end{equation*}
\qed

\begin{Lem} 
\label{dir}
Assume that the commutator $K$ is of trace class, $Kh = \sum_{m = 0}^{\infty}(h, u_m) v_m$, $h \in L^2(\mu)$ where the sum $\sum_{m = 0}^{\infty} \|u_m\| \cdot \|v_m\|$ is finite. Then
$$(W_{+}(\alpha)U - UW_{-}(\alpha))h = \sum_{m = 0}^{\infty} v_m \cdot   \left[(\overline{u_m} h) * \sum_{n = 0}^{\infty} p_{\alpha,n} D_n \right],$$
where $D_n(\zeta) = \sum_{-n}^{n} \zeta^l$ is the Dirichlet kernel of order $n$ and the symbol $*$ denotes the convolution.
\end{Lem}
\beginpf Assume at first that $K = (\cdot, u)v$ and consider the sum $\sum_{l = -n}^{n} U^{l}KU^{-l}$ applied to a vector $h \in L^2(\mu)$:
\begin{multline*}
\sum_{l = -n}^{n} U^{l}KU^{-l} h  = \sum_{l = -n}^{n} U^{l}Kz^{-l}h =  \sum_{l = -n}^{n} U^{l}(z^{-l}h,u)v = \sum_{l = -n}^{n} (z^{-l}h,u)z^{l} v = \\ = \sum_{l = -n}^{n} z^{l} v \int \bar{\xi}^{l}h(\xi) \overline{u(\xi)} \,d\mu(\xi) = v\int h(\xi) \overline{u(\xi)} \sum_{l = -n}^{n}  (z\bar{\xi})^{l} \,d\mu(\xi) = v \cdot [(\bar{u}h) * D_n].
\end{multline*}
By linearity arguments, we obtain the conclusion of the lemma. \qed
\begin{Lem} \label{sym}
Let $k(\xi,z)$ be a real-valued function from $L^\infty(\mu \times \mu)$ such that $k(\xi,z) = k(z,\xi)$. Then the operator $L$ given by the bilinear form 
$$(Lf,g) = \iint (f(\xi) - f(z))\overline{g(z)}k(\xi,z) d\,\mu(\xi)\;d\,\mu(z)$$ 
is a selfadjoint bounded operator.  
\end{Lem} 
\beginpf The quadratic form of the operator $L$ is  
$$(Lf, f) = \iint (f(\xi) - f(z))\overline{f(z)}k(\xi,z)\,d\mu(\xi)\,d\mu(z),
$$ 
and, symmetrically,
$$(Lf, f) = \iint (f(z) - f(\xi))\overline{f(\xi)}k(z,\xi)\,d\mu(z)\,d\mu(\xi).
$$ 
Therefore the sum $2(Lf, f) = - \iint |f(\xi) - f(z)|^2 k(\xi,z)\,d\mu(z)\,d\mu(\xi)$ is real, which proves the statement. \qed

\noindent \textbf{Proof of Theorem \ref{thm1}.} Assume that a triple $\mu$, $U$, $A$ satisfies conditions $1) - 3)$ from Section \ref{introd}. This assumption makes no loss of generality, see \cite[Theorem 7.2]{vk-zns}. 
The theorem will be proved if we check the formula $\lim_{\alpha}((W_{+}(\alpha) - W_{-}(\alpha))h_1,h_2) = 0$ for every pair of vectors $h_1, h_2 \in L^2(\mu)$. A simple computation shows that 
$$U^{-1}W_{+}(\alpha)U = W_{+}(\alpha) - p_{\alpha,0}A + \sum_{n = 0}^{\infty} (p_{\alpha,n} - p_{\alpha,n+1})U^{-n-1}AU^{n+1}.$$ 
By this formula and properties $2)$, $3)$ from the definition of the s-regular summation method, the weak limits of $W_{+}(\alpha) - W_{-}(\alpha)$ and  $W_{+}(\alpha)U - UW_{-}(\alpha)$ do or do not exist simultaneously. Moreover, if they exist, they are or are not equal to zero simultaneously.
By Lemma \ref{dir} we have 
$$
(W_{+}(\alpha)U - UW_{-}(\alpha))h_1 = 
(\phi h_1) * \sum_{n = 0}^{\infty} p_{\alpha,n} D_n -
\phi \cdot   \left[h_1 * \sum_{n = 0}^{\infty} p_{\alpha,n} D_n \right]
$$ 
on every element $h_1 \in L^2(\mu)$. Since the set of all vectors $h_1 \in L^2(\mu)$, that satisfy the condition $\lim_{\alpha}((W_{+}(\alpha)U - UW_{-}(\alpha))h_1,h_2) = 0$ for every element $h_2 \in L^2(\mu)$, forms a reducing subspace of $U$, without loss of generality one can consider $h_1 = 1$. Set $k_{\alpha}(z,\xi) = \sum_{n = 0}^{\infty} p_{\alpha,n} D_n (z \bar{\xi})$. We have  
\begin{equation}
\label{dif1}
p_\alpha(z) = (W_{+}(\alpha)U - UW_{-}(\alpha))1 = \int (\phi(\xi) - \phi(z))k_\alpha(z,\xi)\,d\mu(\xi).
\end{equation}
The norms of $p_\alpha(z)$ in $L^2(\mu)$ are uniformly bounded by $2\|A\|$. Hence one can check the required condition $\lim_{\alpha}(p_\alpha(z),h_2) = 0$, $h \in L^2(\mu)$, only on a dense subset of $L^2(\mu)$. Take the set of all smooth functions on $\T$, which is dense in $L^2(\mu)$. Consider
$$(p_\alpha(z),h_2) = \iint (\phi(\xi) - \phi(z))\overline{h_2(z)}k_\alpha(z,\xi)\,d\mu(\xi)\,d\mu(z),$$
where $h_2 \in C^1(\T)$.  By Lemma \ref{sym} we have
$$(p_\alpha(z),h_2) = \iint \phi(z) (\overline{h_2(\xi)} - \overline{h_2(z)}) k_\alpha(z,\xi)\,d\mu(\xi)\,d\mu(z).$$
Since $|D_n(\bar{\xi} z)| \le \frac{2}{|\xi - z|}$ for every $n \in \N$ (see formula \eqref{dirkernel} below), and $\sum_n p_{\alpha,n} = 1$, the integrand is bounded by $2 |\phi(z)| \sup_{\zeta \in \T}|h^{'}_{2}(\zeta)|$ for every $\alpha \in E$ and $z \in \T$. 

By property $d)$ of s-regular summation method, the averages of every nonconstant unimodular geometric progression tend to zero. For the Dirichlet kernel and $\xi \neq z$ we have 
\begin{equation} \label{dirkernel}
D_n(\bar{\xi}z) = \sum_{-n}^{n} (\bar{\xi}z)^{l} =  \frac{2 \Re((\bar{\xi}z)^{n} - (\bar{\xi}z)^{n+1})}{|1 - \bar{\xi}z|^2}.
\end{equation}
Therefore, the limit $\lim_{\alpha} k_{\alpha}(z,\xi)$ equals zero if $\xi \neq z$. By the dominated convergence theorem we obtain $\lim_{\alpha} (p_\alpha(z),h_2) = 0$.
\qed\\

\textit{Remark.} We have proved the fact that the operators $W_{+}(\alpha) - W_{-}(\alpha)$ tend to zero in the weak operator topology without assuming that the wave operators $W_{\pm}$ exist. 

\section{Proof of Theorem \ref{thm2}.}
\label{proofthm2}
Let $\mu$, $U$, $A$ be triple satisfying conditions $1)-3)$, and let the operators $W_{\pm}(\alpha)$ be defined by formula \eqref{wo}. For the proof of Theorem \ref{thm2} we need the following lemma:
\begin{Lem} \label{conv}
Let $\{h_n\}_{n \in \N}$ be a sequence of elements from $L^2(\mu)$ such that $h_n(z) \to h(z) \in L^2(\mu)$ for $\mu$-almost all points $z \in \T$, and such that the limit $\lim_{\alpha}(W_{-}(\alpha)h_n,h_n)$ exists for all $n \in \N$. Then the limit $\lim_{\alpha}(W_{-}(\alpha)h,h)$ exists. Moreover, if $\sup_{n}\|h_n\| < \infty$, then the double limit $\lim_{\alpha,n}(W_{-}(\alpha)h_n,h_n)$ exists. The same holds for the families $W_{+}(\alpha)$ and $W_{+}(\alpha)  - W_{-}(\alpha)$. 
\end{Lem}
\beginpf In the proof we consider only the case of operators $W_{-}(\alpha)$, other cases can be proven in a similar way. 

Fix an arbitrary $\eps > 0.$ By the Egorov theorem, one can choose a set $F_\eps \subset \T$ such that $\mu(\T \setminus  F_\eps) \le \eps$ and $h_n$ tend to $h$ uniformly on $F_\eps$. Set $h_{n}^{\eps} = \chi_{F_\eps}h_n$ and $h^{\eps} = \chi_{F_\eps}h$, where $\chi_{F_\varepsilon}$ is the characteristic function of the set $F_{\eps}$. Since $h_{n}^{\eps}$ belongs to the reducing subspace of $U$ generated by $h_{n}$, the limits $\lim_{\alpha}(W_{-}(\alpha)h_{n}^{\eps},h_{n}^{\eps})$ exist for every $n \in \N$.
Next, let $n(\eps)$ be the number for which $\|h^{\eps} - h_{n(\eps)}^{\eps}\| < \eps$. Choose the element  $\alpha(\eps) \in E$ such that $|((W_{+}(\alpha_1) - W_{+}(\alpha_2))h_{n(\eps)},h_{n(\eps)})| < \varepsilon$ for every $\alpha_1, \alpha_2 \succ \alpha(\varepsilon)$. We have
\begin{multline*}
|((W_{-}(\alpha_1) - W_{-}(\alpha_2))h,h)| \le  4\|A\|\|h\|\eps +  
|((W_{-}(\alpha_1) - W_{-}(\alpha_2))h^{\eps},h^{\eps})| \le \\
\le 8\|A\|\|h\|\eps +  
|((W_{-}(\alpha_1) - W_{-}(\alpha_2))h_{n(\eps)}^{\eps},h_{n(\eps)}^{\eps}) | \le (1 + 8\|A\|\|h\|)\eps,
\end{multline*}
which implies the existence of $\lim_{\alpha}(W_{+}(\alpha)h,h)$. 

Assume now that norms of $h_n$ in $L^2(\mu)$ are uniformly bounded by $C$. For every elements $\alpha_1, \alpha_2 \succ \alpha(\varepsilon)$ and natural numbers $n_1, n_2 > n(\eps)$ we have 
\begin{multline*}
|(W_{-}(\alpha_1)h_{n_1},h_{n_1}) - (W_{-}(\alpha_2)h_{n_2},h_{n_2})| \le \\ \le  4\|A\|C\eps +  |(W_{-}(\alpha_1)h_{n_1}^{\eps},h_{n_1}^{\eps}) - (W_{-}(\alpha_2)h_{n_2}^{\eps},h_{n_2}^{\eps})| \le \\
\le 8\|A\|C\eps +  
|((W_{-}(\alpha_1) - W_{-}(\alpha_2))h_{n(\eps)}^{\eps},h_{n(\eps)}^{\eps}) | \le \\ \le (1 + 8\|A\|C)\eps,
\end{multline*}
which implies the existence of the double limit $\lim_{\alpha,n}(W_{-}(\alpha)h_n,h_n)$.
\qed
\\

\noindent \textbf{Proof of Theorem \ref{thm2}.} We give the proof for the operators $W_{-}(\alpha)$, the same arguments work for operators $W_{+}(\alpha)$ and $W_{+}(\alpha) - W_{-}(\alpha)$. 

It is shown in \cite[Theorem 7.2]{vk-zns} that Theorem \ref{thm2} holds if we omit item $4)$ in the statement. Hence, our aim is to prove the following fact: If the limit $\lim_{\alpha} W_{-}(\alpha)$ exists for every triple $\mu$, $U$, $A$ satisfying conditions $1)-4)$, then this limit exists for every triple satisfying conditions $1)-3)$.

By the Luzin theorem, we can choose compacts $K_n$ such that $\mu(\T \setminus  K_n) \le \frac{1}{n}$ and $\phi$ coincides with $\phi_n \in C(\T)$ on $K_n$. Set $h_n = \chi_{K_n}$ and $\mu_n = \chi_{K_n}\,d\mu$ (as usual, $\chi_{K_n}$ denotes the characteristic function of the set $K_n$).  For every number $n$ we have $$(W_{-}(\alpha)h_n,h_n) = (M_{h_n}W_{-}(\alpha)M_{h_n}1,1)$$ and the operator $M_{h_n}W_{-}(\alpha)M_{h_n}$ is exactly the past wave operator constructed for the operator $A_n = M_{h_n} A M_{h_n}$. Every triple $\mu_n$, $U$, $A_n$ satisfies properties $1)-4)$, because the commutator $A_nU - UA_n$ equals 
$$(\cdot, h_n\phi)h_n - (\cdot, h_n)h_n\phi_n = (\cdot, \phi_n)1 - (\cdot, 1)\phi_n.$$
An application of Lemma  \ref{conv} ends the proof, due to the fact of pointwise convergence  $h_n \to 1$ on the set $\cup_{n} K_{n}$ of full measure. \qed

\section{Proof of Theorem \ref{thm3}.}\label{proofthm3}
\noindent Define operators $P_r$, $r \in (0,1)$, on $L^2(\mu)$ by the formula
$$
P_r: \phi \mapsto \int_{\T} (\phi(\xi) - \phi(z))\frac{1 - r^2}{|1-r\bar{\xi} z|^2} \,d \mu(\xi).
$$ 
\begin{Lem} \label{lost}
There exist a Borel singular measure $\mu$ on $\T$ without atomic masses and a function $\phi \in C(\T)$ such that the norms of $P_r\phi$  in $L^2(\mu)$ are unbounded.
\end{Lem}

\beginpf
Define a function $\phi$ on $\T$ by 
$$\phi(z)=\begin{cases}
\sqrt[4]{\ARG(z)},&\text{if $\ARG(z) \in [0, \frac{\pi}{2}]$;}\\
0,&\text{if $\ARG(z) \in (0, -\frac{\pi}{2}]$;}\\
g(z),&\text{if $\ARG(z) \in (\frac{\pi}{2}, \frac{3\pi}{2})$.}
\end{cases}$$
Here $g(z)$ is an arbitrary nonnegative function such that $\phi \in C(\T)$. Let $\psi \in C(\T)$ be determined by $\psi(z) = \phi(\bar{z})$. For every measure $\mu$ supported on $I = \{z: \ARG(z) \in [-\frac{\pi}{2}, \frac{\pi}{2}]\}$, we have 
$$(P_r \phi,\psi) = \iint \phi(\xi) \psi(z) \frac{1 - r^2}{|1 - r\bar{\xi}z|^2}\,d\mu(\xi)\,d\mu(z),$$
due to the fact that $\phi\psi = 0$ on $I$. Find a family of arcs $I_k \subset \T$, $k \in \Z \setminus \{0\}$, such that
\begin{itemize}
\item[1)] $a_k = \exp(\sgn(k) 2^{-|k|}i)$ is the center of $I_k$, $ \sgn(k) = \frac{k}{|k|}$, 
\item[2)] $\phi(z) \ge \frac{1}{2}\phi(a_k)$ and $\psi(z) \ge \frac{1}{2}\psi(a_k)$ for every $z \in I_k$,
\item[3)] $\sup_{z \in I_k,\xi \in I_{-k}} |z - \xi| \le 2^{-|k|+2}$.
\end{itemize}
Now define the measure $\mu$ on $I$ by the formula $\mu = \sum_{k \in \Z \setminus \{0\}} k^{-2} \mu_{k}$, where probability measures  $\mu_{k}$ are supported on the arcs $I_k$. As follows from the construction, the measure $\mu$ can be chosen to be singular and free of point masses.  For the measure $\mu$ and the number $r_m = 1- 2^{-m}$, $m \ge 1$, we have
\begin{equation*}
\begin{split}
(P_r \phi,\psi) & \ge \frac{\phi(a_m) \psi(a_{-m})}{4m^4} \iint\limits_{I_m \times I_{-m}} \frac{1 - r_{m}^{2}}{|1 - r_m \bar{\xi} z|^2}\,d\mu(\xi)\,d\mu(z) \\& \ge
\frac{\phi(a_m) \psi(a_{-m})}{4m^4} \inf_{\substack{\xi\in I_{-m}\\z\in I_m}} \frac{1-r_m}{(|z - r_m z| + |r_mz - r_m\xi|)^2} \\& \ge \frac{1}{100 m^{4}}2^m,
\end{split}
\end{equation*}
which tends to infinity when $m$ increases. 
 Thus, we have $\sup_{r}(P_{r}\phi,\psi) = \infty$, therefore $P_{r}\phi$ are cannot be bounded in $L^2(\mu)$.    
\qed\\

\textit{Remark.} In fact, we have proved that $P_{r}\phi$ are unbounded even in the norm of $L^1(\mu)$. This makes the convergence of $P_{r}\phi$ impossible in any reasonable topology.\\

\noindent \textbf{Proof of Theorem \ref{thm3}.} Let $\phi$ be the function constructed in Lemma \ref{lost}. Assume that $AU - UA = (\cdot, \phi)1 - (\cdot, 1)\phi$ and consider the averaged wave operators $W_{\pm}(r)$ with respect to the Abel-Poisson summation method. In our notations this means that $E = ((0,1), \le)$, $p_{r,n} = (1-r)r^n$, and the operators $W_{\pm}(r)$ are defined by formula \eqref{wo} from Section \ref{introd} (now we use the letter $r$ in place of $\alpha$). The Abel-Poisson means of the Dirichlet kernels $D_n(\bar{\xi} z)$ are equal to the Poisson kernel $\frac{1 - r^2}{|1-r\bar{\xi} z|^2}$. Formula \eqref{dif1} gives us
$$
(W_{+}(r)U - UW_{-}(r))1 = P_r\phi.
$$
In particular, we have $\sup_{r}\|P_r\phi\|_{L^2(\mu)} \le 2\|A\|$, which contradicts the conclusion of Lemma \ref{lost}.
\qed

\section{Boundary behaviour of Cauchy-type integrals.}\label{bb}
Lemma \ref{lost} is of special interest because of recent result by V.V.Kapustin concerning the boundary behaviour of Cauchy-type integrals. 
To formulate this result, we introduce several standard objects from the theory of pseudocontinuable functions. The reader can find a more detailed discussion in \cite{VKnew}.
   
Given a singular probability measure $\mu$ on $\T$, define an inner function $\theta$ by the formula 
$$
\Re\left(\frac{1 + \theta(z)}{1 - \theta(z)}\right)  = \int \frac{1 - |z|^2}{|1-\bar{\xi} z|^2}\,d\mu(\xi), \quad |z| < 1.
$$ 
The function $\theta$ generates a family of singular probability measures $\{\sigma_\alpha\}_{\alpha \in \T}$ by   
$$
\Re\left(\frac{\alpha + \theta(z)}{\alpha - \theta(z)}\right)  = \int \frac{1 - |z|^2}{|1-\bar{\xi} z|^2}\,d\sigma_\alpha(\xi), \quad |z| < 1.
$$ 
By the definition, we have $\mu = \sigma_1$. With the inner function $\theta$ we associate the subspace $K_\theta = H^2 \ominus \theta H^2$ of the Hardy space $H^2$ in the unit disk of the complex plane. Functions from $K_\theta$ (also referred to as $\theta$-pseudocontinuable functions) have boundary values $\sigma_\alpha$-almost everywhere for every measure $\sigma_\alpha$, see \cite{polt}. \\

\noindent Define operators $C_r$, $r \in (0,1)$ on $L^2(\mu)$ by the formula
$$
C_r: \phi \mapsto \int_{\T} (\phi(\xi) - \phi(z))\frac{1}{1-r\bar{\xi} z} \,d \mu(\xi).
$$ 

\begin{Thm}[Theorem 1.3 of \cite{VKnew}] \label{alphaconv}
Let a function $h \in K_\theta$ coincide $\sigma_\alpha$-almost everywhere with a continuous function $\phi_\alpha$, where $\alpha \neq 1$. Take the function $\phi \in L^2(\mu)$ such that $h = \phi$ $\sigma_1$-almost everywhere. If $\sigma_1$ has no atomic masses, the family $\{C_r\phi\}$ converges in norm of $L^2(\mu)$ to $\frac{\phi_\alpha - \phi}{\alpha - 1}$.
\end{Thm}
As is shown in \cite{vk-zns}, for the Abel-Poisson summation method the general problem of existence of $W_{\pm}$ in the case of rank-two commutator is equivalent to the convergence of $C_r\phi$ for every function $\phi$ corresponding to a commutator $1)-3)$ from Section \ref{introd}. Theorem \ref{alphaconv} establishes the convergence of $C_r\phi$ for functions $\phi$ with a continuous ``transplantation'' $\phi_{\alpha}$. On the other hand, it follows from Theorem \ref{thm2} that we can check the convergence of $\{C_r\phi\}$ only for continuous functions $\phi$ without loss of generality. By the definition we have $\phi_{1} = \phi$. Unfortunately, the convergence in Theorem \ref{alphaconv} fails if we remove the assumption $\alpha \neq 1$.   
\begin{Prop} \label{divergence}
There exists a Borel singular measure $\mu$ on $\T$ without atomic masses and a continuous function $\phi \in C(\T)$ such that the norms of $C_r\phi$  in $L^2(\mu)$ are unbounded.
\end{Prop}

\beginpf Using the identities 
$$\frac{2}{1 - r\bar{\xi} z} - 1 = \frac{1+r\bar{\xi} z}{1 - r\bar{\xi} z} \quad  \mbox{and} \quad \Re\left(\frac{1+r\bar{\xi} z}{1 - r\bar{\xi} z}\right) = \frac{1 - r^2}{|1 - r\bar{\xi} z|^2},$$ 
one can to show that $\|P_r(\phi)\|_{L^2(\mu)} \le 2 \|C_r(\phi)\|_{L^2(\mu)} + \|\phi\|_{L^2(\mu)}$ for every real-valued function $\phi \in L^2(\mu)$. It remains to apply Lemma \ref{lost}. \qed \\
 
As is easily seen from Proposition \ref{divergence}, the operator $C: f \mapsto \lim_r C_rf$ is not a bounded operator on $L^2(\mu)$. On the other hand, the operator $C$ is well defined on smooth functions. One of natural ways to define it on a wider subset of $L^2(\mu)$ could be the following: at first we define the operator $C$ as the limit $\lim_r C_rf$ on all functions $f \in L^2(\mu)$ for which the limit exists in $L^2(\mu)$, and then take its closure of the graph of $C$. Unfortunately, this way does not work, as the following proposition shows.   
\begin{Prop}
The operator $C$, defined on all functions $f$ for which the limit $\lim_r C_rf$ exists in $L^2(\mu)$, is not a closable operator on $L^2(\mu)$.
\end{Prop}
\beginpf By Theorem \ref{alphaconv}, we have $C\phi = 0$ if $\phi$ is the trace on $\T$ of a continuous function $h \in K_\theta$ (in this case $\phi_\alpha = \phi$ almost everywhere with respect to $\sigma_1$). A well-known result by A.B.Aleksandrov says that continuous functions from $K_\theta$ form a dense subset in $K_\theta$. Since the operator taking a function from $K_\theta$ to its boundary values $\sigma_\alpha-$almost everywhere is a unitary operator from $K_\theta$ to $L^2(\sigma_\alpha)$, see \cite{polt,Clark}, the traces of continuous functions from $K_\theta$ are dense in $L^2(\mu)$. Hence the fact that $C$ a closable operator would imply  $C = 0$. But, obviously, $Cf = -\bar{z} \neq 0$ for the function $f(z) \equiv \bar{z}$. \qed

\enddocument